\documentclass[12pt,a4paper]{article}
\usepackage{amsmath,amsthm,amssymb}
\usepackage{mathrsfs}
\usepackage{latexsym}
\usepackage{lscape}
\usepackage[margin=2.8cm]{geometry}
\usepackage[all]{xy}\UseComputerModernTips

\newtheorem{thm}{Theorem}[section]
\newtheorem{crl}[thm]{Corollary}

\newtheorem{lmm}[thm]{Lemma}

\newtheorem{prp}[thm]{Proposition}

\theoremstyle{definition}
\newtheorem{dfn}[thm]{Definition}
\newtheorem{exa}[thm]{Example}
\theoremstyle{remark}
\newtheorem*{rem}{Remark}

\newcommand{\BR}{{\mathbb R}}
\newcommand{\BC}{{\mathbb C}}
\newcommand{\BZ}{{\mathbb Z}}
\newcommand{\BP}{{\mathbb P}}

\newcommand{\D}{\mathrm{D}}
\newcommand{\Db}{\mathbf{D}^{\mathrm b}}
\newcommand{\Dc}{\mathbf{D}_{\BR \text{-c}}^{\mathrm b}}

\newcommand{\cF}{\mathcal{F}}
\newcommand{\cG}{\mathcal{G}}
\newcommand{\cV}{\mathcal{V}}
\newcommand{\cE}{\mathcal{E}}
\newcommand{\cS}{\mathcal{S}}
\renewcommand{\L}{\mathscr{L}}

\newcommand{\CF}{{\rm CF}}
\newcommand{\reg}{\mathrm{reg}}
\newcommand{\tr}{\mathrm{tr}}
\newcommand{\id}{\mathrm{id}}
\newcommand{\pt}{\mathrm{pt}}
\DeclareMathOperator{\supp}{\operatorname{supp}}
\DeclareMathOperator{\MS}{{\operatorname{SS}}}
\DeclareMathOperator{\Supp}{\operatorname{Supp}}

\DeclareMathOperator{\Ev}{\operatorname{Ev}}
\newcommand{\RG}{R\varGamma}

\DeclareMathOperator{\SHom}{\mathcal{H}om}
\DeclareMathOperator{\Hom}{\mathrm{Hom}}
\newcommand{\e}{\varepsilon}

\newcommand{\tl}[1]{\widetilde{#1}}
\newcommand{\simto}{\overset{\sim}{\longrightarrow}}

\newcommand{\dsum}{\displaystyle \sum}
\newcommand{\dint}{\displaystyle \int}

\renewcommand{\(}{\left(}
\renewcommand{\)}{\right)}

\newcommand{\longhookrightarrow}{\DOTSB\lhook\joinrel\longrightarrow}

\let\Re=\relax
\let\Im=\relax

\DeclareMathOperator{\Re}{\mathrm{Re}}
\DeclareMathOperator{\Im}{\mathrm{Im}}

\makeatletter
  
  \@addtoreset{equation}{section}
\makeatother

\title{Hyperbolic localization via shrinking subbundles\footnote{
\newline Supported by JSPS KAKENHI Grant Number 15J07993 and the Program for
Leading Graduate Schools, MEXT, Japan}}
\author{Yuichi \textsc{Ike}\footnote{Graduate School of Mathematical Sciences, the University
of Tokyo, 3-8-1, Komaba, Meguro-ku, Tokyo, 153-8914, Japan \newline e-mail: \texttt{ike@ms.u-tokyo.ac.jp, yuichi.ike.1990@gmail.com}
}}

\begin{document}

\maketitle
\begin{abstract}
	We study the Lefschetz fixed point formula for constructible sheaves with higher-dimensional fixed point sets.
	We give another proof to the explicit description of Lefschetz cycles in our previous paper.
	For this purpose, we introduce a new notion of shrinking subbundles and describe Lefschetz cycles by using hyperbolic localization with respect to shrinking subbundles.
\end{abstract}

\section{Introduction}
In our previous paper \cite{IMT15}, we gave the explicit description of Lefschetz cycles (see \cite[Theorem~5.10]{IMT15}).
This enabled us to calculate local contributions in the Lefschetz fixed point formula for constructible sheaves in many cases.
However, in \cite{IMT15}, only expanding subbundles were explored and therefore the proof needed some very technical lemmas.
In this paper, we introduce a new notion of \textit{shrinking subbundles} and prove the theorem in a different and straightforward way by using them.
More precisely, the proof in this paper goes as follows.
First, we show that the constructible function of hyperbolic localization with respect to shrinking subbundles is equal to that with respect to expanding subbundles defined in \cite{IMT15} (see Proposition~\ref{prp:compcf}).
Next, we prove that the Lefschetz cycle is equal to the characteristic cycle of the constructible function of hyperbolic localization with respect to shrinking subbundles (see Proposition~\ref{prp:conic}).
In the last section, we explicitly calculate local contributions in the Lefschetz fixed point formula by our method for some examples.

The strategy of the proof in this paper is motivated by the previous works of Goresky-MacPherson~\cite{G-M-1} and Braden~\cite{Braden}.
Indeed, in \cite{G-M-1} the authors showed that the constructible functions with respect to expanding and shrinking subbundles coincide under some assumptions, which are much stronger than ours. 
Note also Braden~\cite{Braden} showed that hyperbolic localization with respect to special expanding and shrinking bundles are isomorphic to each other in the special setting.

\section{The Lefschetz fixed point formula for constructible sheaves and local contributions}
In this paper, we mainly follow the notation of \cite{K-S}.
Let $X$ be a real analytic manifold.
We denote by $\Dc(X)$ the full subcategory of $\Db(X)$ consisting of bounded complexes of $\BC_X$-modules whose cohomology sheaves are $\BR$-constructible (see \cite[Chapter~VIII]{K-S} for the precise definition). 
Let $\phi \colon X \longrightarrow X$ be a real analytic self map of $X$. 
Let $F$ be an object of $\Dc(X)$ and $\Phi$ be a lift of $\phi$, that is, a morphism $\Phi \colon \phi^{-1}F \longrightarrow F$ in $\Dc(X)$.
If the support $\Supp(F)$ of $F$ is compact, $H^j(X;F)$ is a finite-dimensional vector space over $\BC$ for any $j\in \BZ$.
Hence, in this situation, one can define the following number associated with the pair $(F,\Phi)$.

\begin{dfn}
	One defines the \textit{global Lefschetz number} of the pair $(F,\Phi)$ by 
	\begin{align}
		\tr (F,\Phi) :=\dsum_{j \in \BZ} (-1)^j 
		\tr \( H^j(X;F)
		\overset{\Phi}{\longrightarrow} H^j(X;F) \) \in \BC,
	\end{align}
	where the morphism $H^j(X;F) \overset{\Phi}{\longrightarrow} H^j(X;F)$ is induced by 
	\begin{align}
		F \longrightarrow R\phi_*\phi^{-1}F
		\overset{\Phi}{\longrightarrow} R\phi_*F.
	\end{align}
\end{dfn}

Let $M:=\{x \in X \mid \phi(x)=x \}$ be the fixed point set of $\phi$. 
Consider the diagonal embedding $\delta_X \colon X \longhookrightarrow X \times X$ of $X$ and the closed embedding $h:=(\phi, \id_X) \colon X \longhookrightarrow X \times X$ associated with $\phi$. 
Denote by $\Delta_X$ (resp. $\Gamma_{\phi}$) the image of $X$ under $\delta_X$ (resp. $h$). 
Then $M \simeq \Delta_X \cap \Gamma_{\phi}$ and we obtain a chain of morphisms 
\begin{align*}
	R\Hom_{\BC_X}(F,F) 
	\simeq \ & 
	\RG(X;\delta_X^!(F \boxtimes \D F)) \\
	\longrightarrow & 
	\RG_{\Supp(F) \cap \Delta_X}(X \times X; 
	h_*h^{-1}(F \boxtimes \D F)) \\
	\simeq \ & 
	\RG_{\Supp(F) \cap \Delta_X}(X \times X;
	h_*(\phi^{-1}F \otimes \D F)) \\
	\overset{\Phi}{\longrightarrow} & 
	\RG_{\Supp(F) \cap \Delta_X}(X \times X;
	h_*(F \otimes \D F)) \\
	\longrightarrow & 
	\RG_{\Supp(F) \cap \Delta_X}(X \times X; h_*\omega_X) \\
	\simeq \ & \RG_{\Supp(F) \cap M} (X;\omega_X),
\end{align*}
where $\omega_X \simeq \mathrm{or}_X [\dim X]$ is the dualizing complex of $X$ and $\D F \allowbreak = \allowbreak R\SHom_{\BC_X}(F, \omega_X)$ is the Verdier dual of $F$.
Taking the 0-th cohomology, we get a morphism
\begin{align}\label{eq:2-10}
	\Hom_{\Db(X)}(F,F) \longrightarrow H^0_{\Supp(F) \cap M}(X;\omega_X).
\end{align}

\begin{dfn}[\cite{Kashiwara-2}]\label{dfn:2-2} 
	One denotes by $C(F,\Phi)$ the image of $\id_F$ under the morphism~\eqref{eq:2-10} in $H_{\supp(F) \cap M}^0(X;\omega_X)$ and calls it the characteristic class of $(F,\Phi)$.
\end{dfn}

\begin{thm}[\cite{Kashiwara-2}]\label{thm:Lef}
	If $\Supp (F)$ is compact, then 
	\begin{align}
		\tr (F,\Phi)=\dint_X C(F,\Phi).
	\end{align}
	Here $\int_X \colon H^{\dim X}_c(X;\mathrm{or}_X) \longrightarrow \BC$ is the morphism induced by the integral of differential $(\dim X)$-forms with compact support.
\end{thm}

Let $M= \bigsqcup_{i \in I}M_i$ be the decomposition of $M$ into connected components and
\begin{align*}
	H_{\Supp(F) \cap M}^0(X;\omega_X) 
	& =\bigoplus_{i \in I} H_{\Supp(F) \cap M_i}^0(X;\omega_X), \\
	C(F,\Phi) & =\bigoplus_{i \in I} C(F,\Phi)_{M_i}
\end{align*} 
be the associated decomposition.

\begin{dfn}\label{dfn:2-4}
	If $\Supp(F) \cap M_i$ is compact, one sets 
	\begin{align}
		c(F,\Phi)_{M_i}:=\dint_X C(F,\Phi)_{M_i}
	\end{align}
	and calls it the \textit{local contribution} of $(F,\Phi)$ from $M_i$.
\end{dfn}

By Theorem~\ref{thm:Lef}, if $\Supp (F)$ is compact, the global Lefschetz number of $(F,\Phi)$ is the sum of local contributions:
\begin{align}
	\tr (F,\Phi)=\dsum_{i \in I} c(F,\Phi)_{M_i}.
\end{align}
Hence, the next important problem is to calculate the local contribution $c(F,\Phi)_{M_i}$.

\section{Lefschetz cycles and microlocal index formula}\label{sec:def}

In this section, we recall Lefschetz cycles defined in Matsui-Takeuchi~\cite{M-T-3}. 
Assume that the fixed point set $M=\{ x \in X \mid \phi(x)=x\}$ of $\phi \colon X \longrightarrow X$ is a smooth submanifold of $X$.
Then the differential of $\phi$ induces the endomorphism $\phi'$ of the normal bundle $T_MX$: $\phi' \colon T_MX \longrightarrow T_MX$.

\begin{dfn}
	Let $V$ be a finite-dimensional vector space over $\BR$.
	For an $\BR$-linear endomorphism $f \colon V \longrightarrow V$, one sets
	\begin{align}
		\Ev(f):=\{\text{the eigenvalues of } f^\BC \colon
		V^\BC \longrightarrow V^\BC \} \subset \BC,
	\end{align}
	where $V^\BC$ is the complexification of $V$.
\end{dfn}

We also assume that
\begin{align}
	\text{``$1\notin \Ev(\phi^{\prime}_x)$ for any $x \in M$"}.
\end{align}
Note that this assumption is equivalent to 
\begin{align*}
	\text{``$\Gamma_{\phi}=\{ (\phi(x),x) \in X \times X \mid x \in X\}$ intersects with $\Delta_X$ cleanly along $M$ in $X \times X$"}.
\end{align*}
Under the above assumptions, one can show that $\cF=T^*_{\Gamma_{\phi}}(X\times X) \cap T^*_{\Delta_X}(X\times X)$ is a vector bundle over $M \simeq \Gamma_\phi \cap \Delta_X$ and isomorphic to $T^*M$.
In this situation, Matsui-Takeuchi~\cite{M-T-3} constructed a Lagrangian cycle $LC(F,\Phi) \in H^0_{\MS(F) \cap \cF}(\cF;\pi_M^{-1}\omega_M)$ in $\cF \simeq T^*M$, called the \textit{Lefschetz cycle} of $(F,\Phi)$.
Here $\MS(F) \subset T^*X$ denotes the micro-support of $F$.
One can calculate local contributions of $(F,\Phi)$ as intersection numbers of continuous sections of $\cF \simeq T^*M$ and $LC(F,\Phi)$ as follows.
Let $M =\bigsqcup_{i \in I}M_i$ be the decomposition of $M$ into connected components.
Then we get a decomposition $\cF \simeq \bigsqcup_{i \in I} \cF_i$ with $\cF_i\simeq T^*M_i$.
Moreover we obtain the associated decomposition
\begin{align}
	LC(F,\Phi) =\dsum_{i \in I} LC(F,\Phi)_{M_i}
\end{align}
of $LC(F,\Phi)$, where 
$LC(F,\Phi)_{M_i}$ 
is a Lagrangian cycle in $\cF_i \simeq T^*M_i$. 
Now let us fix a fixed point component $M_i$ and assume that $\Supp(F) \cap M_i$ is compact. 
We shall show how the local contribution $c(F,\Phi)_{M_i}\in \BC$ of $(F,\Phi)$ from $M_i$ can be expressed by $LC(F,\Phi)_{M_i}$. 
For simplicity, we denote $M_i$, $\cF_i$, $LC(F,\Phi)_{M_i}$, $c(F,\Phi)_{M_i}$ by $M$, $\cF$, $LC(F,\Phi)$, $c(F,\Phi)$, respectively. 
With any continuous section $\sigma \colon M \longrightarrow \cF \simeq T^*M$ of the vector bundle $\cF$, we can associate a cycle $[\sigma] \in H^0_{\sigma(M)}(\cF;\pi_M^!\BC_M)$ which is the image of $1 \in H^0(M;\BC_M)$ by the isomorphism $H^0_{\sigma(M)}(T^*M;\pi_M^!\BC_M) \simeq H^0(M;(\pi_M\circ \sigma)^!\BC_M) \simeq H^0(M;\BC_M)$ (see \cite[Definition~9.3.5]{K-S}).
If $\sigma(M) \cap \supp (LC(F,\Phi))$ is compact, we can define the intersection number $\# ([\sigma] \cap LC(F,\Phi))$ of $[\sigma]$ and $LC(F,\Phi)$ to be the image of $[\sigma] \otimes LC(F,\Phi)$ by the chain of morphisms 
\begin{align*}
	H^0_{\sigma(M)}(\cF;\pi_M^!\BC_M) \otimes H^0_{\supp(LC(F,\Phi))}(\cF;\pi_M^{-1}\omega_M) 
	\longrightarrow \ & H^0_{\sigma(M) \cap \supp (LC(F,\Phi))}(\cF;\omega_{\cF}) \\
	\overset{\int_{\cF}}{\longrightarrow} \ & \BC.
\end{align*}

\begin{thm}[{\cite[Theorem~4.8]{M-T-3}}]\label{thm:idx}
	Assume that $\Supp (F) \cap M$ is compact. 
	Then for any continuous section $\sigma \colon M \longrightarrow \cF \simeq T^*M$ of $\cF$, one has an equality
	\begin{align}
		c(F,\Phi) =\# ([\sigma] \cap LC(F,\Phi)).
	\end{align}
\end{thm}

\begin{rem}
Guillermou~\cite{Guillermou} defined the microlocal Lefschetz classes of elliptic pairs in more general situations.
Note also that if $\phi=\id_X$, $M=X$ and $\Phi=\id_F$, the Lagrangian cycle $LC(F,\Phi)$ coincides with the characteristic cycle $CC(F)$ of $F$ introduced by Kashiwara \cite{Kashiwara-1}. 
For recent results on this subject, see also \cite{K-S-2}, \cite{Ike}, etc.
\end{rem}

As a basic property of Lefschetz cycles, we have the following homotopy invariance. 
Let $I=[0,1]$ and let $\phi \colon X\times I \longrightarrow X$ be the restriction of a morphism of real analytic manifolds $X \times \BR \longrightarrow X$. 
For $t \in I$, let $i_t \colon X \longhookrightarrow X \times I$ be the map defined by $x \longmapsto (x,t)$ and set $\phi_t :=\phi \circ i_t \colon X \longrightarrow X$. 
Assume that the fixed point set of $\phi_t$ in $X$ is a smooth submanifold of $X$ and does not depend on $t \in I$. 
We denote this fixed point set by $M$. 
Let $F\in \Dc(X)$ and $\Phi \colon \phi^{-1}F \longrightarrow p^{-1}F$ be a morphism in $\Dc(X \times I)$, where $p
\colon X \times I \longrightarrow X$ is the projection. 
We set
\begin{align}
	\Phi_t :=\Phi|_{X\times \{t\}} \colon \phi_t^{-1}F \longrightarrow F
\end{align}
and 
\begin{align}
	\cF_t:=T^*_{\Gamma_{\phi_t}}(X\times X) \cap T^*_{\Delta_X}(X\times X) \simeq T^*M
\end{align}
for $t \in I$. 

\begin{prp}\label{prp:ht}
	Assume that $\Supp(F) \cap M$ is compact and $\MS(F) \cap \cF_t \subset T^*M$ does not depend on $t \in I$ as a subset of $T^*M$. 
	Then the Lefschetz cycle $LC(F,\Phi_t)$ does not depend on $t \in I$.
\end{prp}

Next, we shall give a formula which enables us to describe the Lefschetz cycle $LC(F,\Phi)$ explicitly in the special case where $\phi \colon X \longrightarrow X$ is the identity map of $X$ and $M=X$.
For this purpose, until the end of this section, we consider the situation where $\phi =\id_X$, $M=X$ and $\Phi \colon F \longrightarrow F$ is an endomorphism of $F \in \Dc(X)$. 
In this case, $LC(F,\Phi)$ is a Lagrangian cycle in $T^*X$. 
For a real analytic function $f \colon Y \longrightarrow \BR$ on a real analytic manifold $Y$, we define a section $\sigma_{f} \colon Y \longrightarrow T^*Y$ of $T^*Y$ by $\sigma_{f}(y):=(y ;df(y))\ (y \in Y)$ and set 
\begin{align}
	\Lambda_{f} :=\sigma_{f}(Y)=\{ (y ;df(y)) \mid y \in Y \}.
\end{align}
Note that $\Lambda_{f}$ is a Lagrangian submanifold of $T^*Y$.

\begin{thm}\label{thm:localinter}
	Let $X$, $F \in \Dc(X)$ and $\Phi \colon F \longrightarrow F$ be as above.
	For a real analytic function $f \colon X \longrightarrow \BR$ and a point $x_0 \in X$, assume that
	\begin{align}
		\Lambda_{f} \cap \MS(F) \subset \{(x_0;df(x_0))\}.
	\end{align}
	Then the intersection number $\# ([\sigma_{f}] \cap LC(F,\Phi))$ (at the point $(x_0;df(x_0))\in T^*X$) is equal to
	\begin{align}
		\dsum_{j\in \BZ}(-1)^j \tr \( H^j_{\{f \geq f(x_0)\}}(F)_{x_0} \overset{\Phi}{\longrightarrow}H^j_{\{f \geq f(x_0)\}}(F)_{x_0} \).
	\end{align}
\end{thm}

By Theorem~\ref{thm:localinter}, we can describe the Lefschetz cycle $LC(F,\Phi)$ as follows. 
Let $X=\bigsqcup_{\alpha \in A}X_{\alpha}$ be a $\mu$-stratification of $X$ such that
\begin{align}
	\supp(LC(F,\Phi)) \subset \MS(F) \subset \bigsqcup_{\alpha \in A} T_{X_{\alpha}}^*X.
\end{align}
Then $\Lambda := \bigsqcup_{\alpha \in A} T_{X_{\alpha}}^*X$ is a closed conic subanalytic Lagrangian subset of $T^*X$. Moreover there exists an open dense smooth subanalytic subset $\Lambda_0$ of $\Lambda$ whose decomposition $\Lambda_0=\bigsqcup_{i \in I}\Lambda_i$ into connected components satisfies the condition
\begin{align}
	\text{``for any $i\in I$, there exists $\alpha_i \in A$ such that $\Lambda_i
\subset T_{X_{\alpha_i}}^*X$"}.
\end{align}

\begin{dfn}
	For $i \in I$ and $\alpha_i\in A$ as above, one defines a complex number $m_i \in \BC$ as follows.
	Take a point $p \in \Lambda_i$ and set $x:=\pi_X(p) \in \pi_X(\Lambda_i) \subset X_{\alpha_i}$. 
	Take moreover a real analytic function $f \colon X \longrightarrow \BR$ defined in an open neighborhood of $x$ in $X$ which satisfies the following conditions:
	\begin{enumerate}
		\item $p=(x;df(x)) \in \Lambda_i$,
		\item the Hessian $\mathrm{Hess} (f|_{X_{\alpha_i}})$ of $f|_{X_{\alpha_i}}$ is positive definite.
	\end{enumerate}
	Then one sets
	\begin{align}
		m_i 
		:=
		\dsum_{j \in \BZ} (-1)^j \tr\(H^j_{\{f \geq f(x)\}}(F)_x \overset{\Phi}{\longrightarrow} H^j_{\{f \geq f(x)\}}(F)_x\).
	\end{align}	
\end{dfn}

\begin{crl}\label{crl:coeff}
	In the situation as above, for any $i \in I$ one has 
	\begin{align}
		LC(F,\Phi)=m_i \cdot [T_{X_{\alpha_i}}^*X]
	\end{align}
	in an open neighborhood of $\Lambda_i$ in $T^*X$.
\end{crl}

Now let us define a $\BC$-valued constructible function $\varphi(F,\Phi)$ on $X$ by
\begin{align}
	\varphi(F,\Phi)(x)
	:=\dsum_{j \in \BZ} (-1)^j \tr \( H^j(F)_x \xrightarrow{\Phi|_{\{x\}}} H^j(F)_x\) \label{eq:defphi}
\end{align}
for $x\in X$. 
Recall that for a real analytic manifold $Z$ there exists an isomorphism called the characteristic cycle map 
\begin{align}
	CC \colon \CF(Z)_{\BC} \simto \varGamma(T^*Z;\L_Z)
\end{align}
between the $\BC$-vector space consisting of the $\BC$-valued constructible functions on $Z$ and that of closed conic subanalytic Lagrangian cycles on $T^*Z$ with coefficients in $\BC$.
By Corollary~\ref{crl:coeff}, we have the following.

\begin{thm}\label{thm:5-8}
	In the situation $\phi=\id_X$, $\Phi \colon F \longrightarrow F$, etc.\ as above, one has an equality
	\begin{equation}
		LC(F,\Phi)=CC(\varphi(F,\Phi))
	\end{equation}
	as Lagrangian cycles in $T^*X$.
\end{thm}

\section{Hyperbolic localization by using shrinking subbundles}

In this section, we explicitly describe the Lefschetz cycle $LC(F,\Phi)$ in many cases. 
Let $M$ be a fixed point component of $\phi \colon X \longrightarrow X$ and assume that $M$ is a smooth submanifold of $X$ for simplicity (actually all we need to assume is $\Supp(F) \cap M \subset M_{\reg}$). 
Throughout this section, we assume the condition
\begin{align}
	\text{``$1\notin \Ev(\phi^{\prime}_x)$ for any $x \in M$"}.
\end{align}
Here, recall that $\phi' \colon T_MX \longrightarrow T_MX$ is the endomorphism of the normal bundle $T_MX$ induced by $\phi$.
Under the assumptions, $\cF=T^*_{\Gamma_{\phi}}(X\times X) \cap T^*_{\Delta_X}(X\times X)$ is a vector bundle over $M$ which is isomorphic to $T^*M$, and we can define the Lefschetz cycle $LC(F,\Phi)$ in $\cF$. 
Recall also that the characteristic cycle map induces an isomorphism
\begin{align}
	CC \colon \CF(M)_{\BC} \simto \varGamma(T^*M;\L_M).
\end{align}
Thus there exists a unique constructible function $\theta_M \in \CF(M)_{\BC}$ on $M$ such that 
\begin{align}
	LC(F,\Phi)=CC(\theta_M).
\end{align}
The aim of this section is to construct the function $\theta_M$ explicitly.

Under the above assumptions, the fixed point set of $\phi^{\prime} \colon T_MX \longrightarrow T_MX$ is the zero-section $M$. 
Let $\Gamma_{\phi^{\prime}}:= \{(\phi^{\prime}(p),p) \mid p \in T_MX\}\subset T_MX \times T_MX$ be the graph of $\phi^{\prime}$ and $\Delta_{T_MX}\simeq T_MX$ the diagonal subset of $T_MX \times T_MX$. 
Then
\begin{align}
	\cF^{\prime}:=T^*_{\Gamma_{\phi^{\prime}}}(T_MX\times T_MX) \cap T^*_{\Delta_{T_MX}}(T_MX \times T_MX)
\end{align}
is a vector bundle over the zero-section $M \simeq \Gamma_{\phi^{\prime}} \cap \Delta_{T_MX}$ of $T_MX$. 
Since $\cF^{\prime}$ is also isomorphic to $T^*M$ by our assumptions, we shall identify it with 
$\cF=T^*_{\Gamma_{\phi}}(X \times X) \cap T^*_{\Delta_X}(X\times X)$.
Now consider the specialization $\nu_M(F) \in \Dc(T_MX)$ of $F$ along $M$ (cf.\ \cite[Section 4.2]{K-S}).
Note that $\nu_M(F)$ is a conic object on $T_MX$.
There is a natural morphism
\begin{align}
	\Phi^{\prime} \colon (\phi^{\prime})^{-1}\nu_M(F) \longrightarrow \nu_M(F)
\end{align}
induced by $\Phi \colon \phi^{-1}F \longrightarrow F$. 
Hence from the pair $(\nu_M(F), \Phi^{\prime})$, we can construct the Lefschetz cycle $LC(\nu_M(F),\Phi^{\prime})$ in $\cF' \simeq \cF$.

\begin{prp}[{\cite[Proposition~5.1]{IMT15}}]\label{prp:specialization}
	In $\cF \simeq \cF^{\prime}$, one has
	\begin{align}
		LC(F,\Phi) =LC(\nu_M(F),\Phi^{\prime}).
	\end{align}
\end{prp}

From now on, we shall identify $\cF \simeq \cF^{\prime}$ with $T^*M$ and describe $LC(F,\Phi)=LC(\nu_M(F),\Phi^{\prime})$.
Since our result holds for any conic object on any vector bundle over $M$, we consider the following general setting. 
Let $\tau \colon \cG \longrightarrow M$ be a real vector bundle of rank $r>0$ over $M$ and $\psi \colon \cG \longrightarrow \cG$ its endomorphism. 
Assume that the fixed point set of $\psi$ is the zero-section $M$ of $\cG$. 
This assumption implies that 
\begin{align}
	\text{``$1 \notin \Ev(\psi_x)$ for any $x \in M$"}. \label{cd:6-1}
\end{align}
Suppose that we are given a conic $\BR$-constructible object $G \in \Dc(\cG)$ on $\cG$ and a morphism $\Psi \colon \psi^{-1}G \longrightarrow G$ in $\Dc(\cG)$. 
Since the fixed point set of $\psi$ is the zero-section $M$, we have an isomorphism
\begin{align}
	\cF_0:=T^*_{\Gamma_{\psi}}(\cG \times \cG) \cap T^*_{\Delta_{\cG}}(\cG \times \cG) \simeq T^*M.
\end{align}
Note that the Lefschetz cycle $LC(G,\Psi)$ is a Lagrangian cycle in $\cF_0 \simeq T^*M$.
\smallskip

Fix a point $\mathring{x} \in M$ and consider the linear homomorphism $\psi_{\mathring{x}} \colon \cG_{\mathring{x}}
 \longrightarrow \cG_{\mathring{x}}$.
Let $\lambda_1,\dots ,\lambda_d$ be the eigenvalues of $\psi_{\mathring{x}}$ on $[0,1]$ and $\lambda_{d+1},\dots ,\lambda_r$ the remaining ones. 
Here, the eigenvalues are counted with multiplicity.
Since these eigenvalues vary depending on $x \in M$ continuously, we denote their continuous extensions to a 
neighborhood of $\mathring{x}$ in $M$ by $\lambda_1(x),\dots \lambda_r(x)$. 
Then there exists a sufficiently small $\e>0$ such that 
\begin{align}
	\lambda_{d+1},\dots,\lambda_r \not\in \{ z \in \BC \mid |z| \le 1, \Re z \ge -\e, |\Im z| \le \e \}.
\end{align}
By the continuity of the eigenvalues, there exists a sufficiently small neighborhood $U$ of $\mathring{x}$ in $M$ such that 
\begin{align}
	\lambda_{d+1}(x),\dots, \lambda_r(x) & \not\in \{ z \in \BC \mid |z| \le 1, \Re z \ge -\e, |\Im z| \le \e \}, \\
	\lambda_1(x),\dots, \lambda_d(x) & \in \{ z \in \BC \mid |z| \le 1, \Re z \ge -\e, |\Im z| \le \e \} 
\end{align}
for any $x \in U$.
If necessary, replacing $U$ by a smaller one, we may assume also that $\cG$ is trivial on $U$. 
For $x \in U$ we set 
\begin{align}
	P_x= \frac{1}{2 \pi i} \int_{C}(z- \psi_x)^{-1} dz,\label{eq:projector}
\end{align}
where $C$ is the path on the boundary of the set $\{ z \in \BC \mid |z| \le 1, \Re z \ge -\e, |\Im z| \le \e \} \subset \BC$. 
Then $P_x \colon \cG_x \longrightarrow \cG_x$ is the projector onto the direct sum of the generalized eigenspaces associated with the eigenvalues in the set. 
The family $\{P_x \}_{x \in U}$ defines an endomorphism $P$ of $\cG|_U$, whose image $\cV \subset \cG|_U$ is a subbundle of $\cG|_U$.

\begin{dfn}\label{dfn:6-1}
	The subbundle $\cV:=\Im P \subset \cG|_U$ is said to be the \textit{minimal shrinking subbundle}  of $\cG|_U$ (on the neighborhood $U$ of $\mathring{x} \in M$).
\end{dfn}

Recall that for a real vector space $V$, $V^\BC:=V \otimes_{\BR} \BC$ denotes its complexification.
Moreover for an endomorphism $f \colon V \longrightarrow V$ of $V$, we denote by $V^\BC_\lambda$ the generalized eigenspace of $V^\BC$ with generalized eigenvalue $\lambda$ for $\lambda \in \BC$.

\begin{dfn}[{cf.\ \cite[Section 9.6]{K-S}}]\label{dfn:6-2} 
	A subbundle $\cS$ of $\cG|_U(=:\tl{U})$ is said to be a \textit{shrinking subbundle} if it satisfies the following conditions:
	\begin{enumerate}
		\item $\psi|_{\tl{U}}(\cS) \subset \cS$,
		\item $\cV$ is a subbundle of $\cS$,
		\item $\displaystyle \cS_x^\BC \subset \bigoplus_{\lambda \not\in [1,+\infty)}(\cG_x)_\lambda^\BC$ for any $x \in U$.
	\end{enumerate}
\end{dfn}

Note that for any $\mathring{x} \in M$, there is a sufficiently small open neighborhood $U$ and a shrinking subbundle $\cS$ of $\cG|_U$.
Note also that for any shrinking subbundle $\cS$ of $\cG|_U$, we have $(\psi|_{\tl{U}})^{-1}(\cS)=\cS$.

\begin{dfn}[\cite{Braden}, \cite{K-S}]\label{dfn:6-3} 
	Let $\tau_{\cS} \colon \cS \longrightarrow U$ be a shrinking subbundle of $\cG|_U(=\tl{U})$.
	Consider the sequence of morphisms $U \overset{i_\cS}{\longhookrightarrow} \cS \overset{j_\cS}{\longhookrightarrow} \cG$, where $i_\cS$ is the zero-section of $\cS$ and $j_\cS$ is the inclusion.
	One defines an object $G_\cS^{-1!} \in \Dc(U)$ by 
	\begin{align}
		G_\cS^{-1!}:=i_\cS^{-1} j_\cS^! G \simeq R{\tau_\cS}_*\RG_\cS(G|_{\tl{U}})
	\end{align}
	and its endomorphism $\Psi_\cS^{-1!} \colon G_\cS^{-1!} \longrightarrow G_\cS^{-1!}$ by the composite of the morphisms 
	\begin{align*}
		R{\tau_\cS}_*\RG_\cS(G|_{\tl{U}})
		\longrightarrow & 
		R{\tau_\cS}_*\RG_\cS \tl{\psi}_*\tl{\psi}^{-1}(G|_{\tl{U}}) \\
		\simeq \ & 
		R{\tau_\cS}_*\tl{\psi}_*\RG_\cS \tl{\psi}^{-1}(G|_{\tl{U}})
		\simeq R{\tau_\cS}_* \RG_\cS((\psi^{-1}G)|_{\tl{U}}) \\
		\overset{\Psi}{\longrightarrow} &
		R{\tau_\cS}_*\RG_\cS(G|_{\tl{U}}).
	\end{align*}
	Here one sets $\tl{\psi}:=\psi|_{\tl{U}}$ and the first morphism above is induced by the adjunction.
	The pair $(G_\cS^{-1!},\Psi_\cS^{-1!})$ is called the \textit{hyperbolic localization} of $(G,\Psi)$ with respect to the shrinking subbundle $\cS$.
\end{dfn}

\begin{prp}\label{prp:shrinking}
	Let $\mathring{x} \in M$ be a point of $M$. 
	Then there exists a sufficiently small open neighborhood $U$ of $\mathring{x}$ in $M$ such that for any subanalytic relatively compact open subset $V$ of $U$ and for any shrinking subbundle $\cS$ of $\cG|_U$, one has
	\begin{align}
		\int_{\tl{U}} C(\RG_{\tl V}(G)|_{\tl{U}},\RG_{\tl V}(\Psi)|_{\tl{U}})
		=
		\tr((\RG_{\tl V}(G))_\cS^{-1!},(\RG_{\tl V}(\Psi))_\cS^{-1!}).
	\end{align}
	Here one sets $\tl{U}:=\tau^{-1}(U)$ and $\tl{V}:=\tau^{-1}(V)$.
\end{prp}

The proof of this proposition is similar to that of \cite[Proposition~9.6.12]{K-S} and we omit it here. 
By the isomorphism $(\RG_{\tl V}(G))_\cS^{-1!} \simeq \RG_{V}(G_\cS^{-1!})$, we obtain an equality
\begin{align}\label{ADDEE} 
	\int_{\tl{U}} C(\RG_{\tl V}(G)|_{\tl{U}},\RG_{\tl V}(\Psi)|_{\tl{U}}) 
	=
	\tr(\RG_{V}(G_\cS^{-1!}),\RG_{V}(\Psi_\cS^{-1!})).
\end{align}

Take a sufficiently small open subset $U$ of $M$ for which Proposition~\ref{prp:shrinking} holds and define a constructible function $\varphi(G_\cS^{-1!},\Psi_\cS^{-1!})$ on it associated with the hyperbolic localization $(G_\cS^{-1!},\Psi_\cS^{-1!})$ by 
\begin{align}
	\varphi(G_\cS^{-1!},\Psi_\cS^{-1!})(x)
	:=
	\dsum_{j \in \BZ} (-1)^j\tr\( H^j (G_\cS^{-1!})_x \xrightarrow{\Psi_\cS^{-1!}|_{\{x\}}} H^j(G_\cS^{-1!})_x\). 
\end{align}
For the notation $\varphi(F,\Phi)$, see \eqref{eq:defphi}.
By the constructibility, the value is equal to the alternating sum of the traces on a sufficiently small open ball centered at $x$.
Then by applying \eqref{ADDEE}, we find that it does not depend on the choice of the shrinking subbundle $\cS$. 
Hence we can glue such locally defined constructible functions to obtain a global one $\varphi^s_M(G,\Psi)$ on $M$.
\smallskip

Next, we shall compare the two constructible functions associated with expanding subbundles and shrinking subbundles.
In \cite{IMT15}, we defined the constructible function associated with expanding subbundles $\varphi^e_M(G,\Psi)$ (in \cite{IMT15} we used the symbol $\varphi_M(G,\Psi)$).
Recall that a subbundle $\cE$ of $\cG|_U$ on a sufficiently small open set $U$ is said to be an expanding subbundle 
if it satisfies the following conditions:
\begin{enumerate}
	\item $\psi|_{\tl{U}}(\cE) \subset \cE$.
	\item The minimal expanding subbundle is contained in $\cE$ as a subbundle.
	\item $\displaystyle \cE_x^\BC \subset \bigoplus_{\lambda \not\in [0,1]} (\cG_x)_\lambda^\BC$ for any $x \in U$.
\end{enumerate}
Let $\tau_\cE \colon \cE \longrightarrow U$ be an expanding subbundle of $\cG|_U$ and denote by $i_{\cE} \colon U \longhookrightarrow \cE$ its zero-section.
The hyperbolic localization of $(G,\Psi)$ with respect to the expanding subbundle $\cE$ is the pair $(G_\cE^{!-1},\Psi_\cE^{!-1})$ of 
\begin{align}
	G_\cE^{!-1}:=i_\cE^!(G|_\cE) \simeq R{\tau_\cE}_!(G|_\cE) \ \in \Dc(U)
\end{align}
and the associated endomorphism $\Psi_\cE^{!-1} \colon G_\cE^{!-1} \longrightarrow G_\cE^{!-1}$.
See \cite[Section 5]{IMT15} for the precise definition.
The value of the function $\varphi^e_M(G,\Psi)$ at $x \in M$ is defined as follows.
Taking an expanding subbundle on a sufficiently small neighborhood of $x$, we set 
\begin{align}
	\varphi^e_M(G,\Psi)(x)
	:=\dsum_{j \in \BZ} (-1)^j\tr\( H^j (G_\cE^{!-1})_x \xrightarrow{\Psi_\cE^{!-1}|_{\{x\}}} H^j(G_\cE^{!-1})_x\). 
\end{align}

\begin{prp}\label{prp:compcf}
	One has an equality 
	\begin{align}
		\varphi_M^e(G,\Psi)
		=
		\varphi_M^s(G,\Psi)
	\end{align}
	as elements of $\CF(M)_\BC$.
\end{prp}

\begin{proof}
	Fix $x_0 \in M$ and we shall compare the values at $x_0$.
	Taking a sufficiently small neighborhood $U$ of $x_0$, we may assume that there exist the minimal expanding subbundle and the minimal shrinking subbundle on $U$. 
	Moreover, by the homotopy invariance of traces,  replacing $\psi$ by $t \psi$ with $|1-t| \ll 1$ and shrinking $U$ if it is necessary, we may assume that
	\begin{align}
		\text{``$\Ev(\psi_x) \cap \{z \in \BC \mid |z|=1\}=\emptyset$ for any $x \in U$"}.
	\end{align}
	from the first.
	By taking the unit circle $\{z \in \BC \mid |z|=1\}$ as the integration path $C$ in \eqref{eq:projector}, we can construct the projector onto the direct sum of the generalized eigenspaces associated with the eigenvalues in the unit ball. 
	Using the projector, we can decompose the vector bundle $\cG|_U$ on $U$ as $\cG|_U \simeq \cG_+ \oplus \cG_-$ satisfying
	\begin{align}
		(\cG_{+})_x
		=
		\left(\bigoplus_{|\lambda|>1}(\cG_{x})_\lambda^\BC \right)\cap \cG_x, 
		\quad
		(\cG_{-})_x
		=
		\left(\bigoplus_{|\lambda|<1}(\cG_{x})_\lambda^\BC \right) \cap \cG_x.
	\end{align}
	Furthermore, by shrinking $U$ if it is necessary, we endow $\cG|_U$ with a metric such that 
	\begin{equation}\label{cd:metric}
		\begin{split}
			& \text{``there exist constants $c_1, c_2$ with $0<c_1<1<c_2$ satisfying the condition} \\
			& \text{$|\psi_{x}(v_-)| \le c_1|v_-| \ (v_- \in \cG_{-,x}), \ | \psi_{x}(v_+)| \ge c_2|v_+| \ (v_+ \in \cG_{+,x})$ for any $x \in U$"}. 
		\end{split}
	\end{equation}
	In what follows, we write $G$ for $G|_{\tl{U}}$ and $\Psi$ for $\Psi|_{\tl{U}}$, for simplicity. 
	Using the metric, we define the subset $Z_{a,b} \subset \cG|_U$ by 
	\begin{align}
		Z_{a,b}
		:=
		\{
			(x,v_+,v_-) \in \cG_+ \oplus \cG_- \mid |v_+|<a, |v_-| \le b 
		\}
	\end{align}
	for some constants $a, b >0$.
	In what follows, we fix $a,b>0$ and write $Z$ instead of $Z_{a,b}$.

	By \cite[Proposition~5.5]{IMT15} (or \cite[Proposition~9.6.12]{K-S}), the value of $\varphi_M^e(G,\Psi)$ at $x_0$ is calculated as 
	\begin{align}
		\varphi_M^e(G,\Psi)(x_0)
		=
		\int_{\cG_x} C(G|_{\tau^{-1}(x_0)},\Psi|_{\tau^{-1}(x_0)}).
	\end{align}
	By the local invariance of the characteristic class and the Lefschetz theorem (Theorem~\ref{thm:Lef}), we have
	\begin{align*}
		\int_{\cG_x} C(G|_{\tau^{-1}(x_0)},\Psi|_{\tau^{-1}(x_0)}) 
		& =
		\int_{\cG_x} C(G_Z|_{\tau^{-1}(x_0)},\Psi_Z|_{\tau^{-1}(x_0)}) \\
		& =
		\tr(G_Z|_{\tau^{-1}(x_0)}, \Psi_Z|_{\tau^{-1}(x_0)}). 
	\end{align*}
	Here, the last equality follows from the fact that $G_Z|_{\tau^{-1}(x)}$ has a compact support.
	We thus obtain an equality
	\begin{align}
		\varphi_M^e(G,\Psi)(x_0)
		=
		\tr(G_Z|_{\tau^{-1}(x_0)}, \Psi_Z|_{\tau^{-1}(x_0)}). \label{eq:cfex}
	\end{align}
	
	On the other hand, by Proposition~\ref{prp:shrinking}, $\varphi_M^s(G,\Psi)(x_0)$ can be calculated as 
	\begin{align}
		\varphi_M^s(G,\Psi)(x_0)
		=
		\int_{\tl{U}} C(\RG_{\tl{V}}(G),\RG_{\tl{V}}(\Psi)),
	\end{align}
	where $V:=B(x_0,\delta)$ is a sufficiently small open ball centered at $x_0$ and $\tl{V}:=\tau^{-1}(V) \subset \tl{U}$.
	By the local invariance and the Lefschetz theorem (Theorem~\ref{thm:Lef}), this is equal to 
	\begin{align}
		\int_{\tl{U}} C(\RG_{\tl{V}}(G_Z),\RG_{\tl{V}}(\Psi_Z)) 
		=
		\tr(\RG_{\tl{V}}(G_Z),\RG_{\tl{V}}(\Psi_Z)).
	\end{align}
	Since $\tau \colon \tl{U} \longrightarrow U$ is proper on $\Supp(G_Z)$, for a sufficiently small $\delta >0$ we have isomorphisms
	\begin{align*}
		\RG(\tl{U};\RG_{\tl{V}}(G_Z)) 
		& \simeq 
		\RG(\tau^{-1}(B(x_0,\delta));G_Z) \\
		& \simeq 
		\RG(\tau^{-1}(x_0);G_Z|_{\tau^{-1}(x_0)}).
	\end{align*}
	Thus, for a sufficiently small $\delta >0$, we obtain
	\begin{align}
		\tr(\RG_{\tl{V}}(G_Z),\RG_{\tl{V}}(\Psi_Z))
		=
		\tr(G_Z|_{\tau^{-1}(x_0)},\Psi_Z|_{\tau^{-1}(x_0)})
	\end{align}
	and 
	\begin{align}
		\varphi^s_M(x_0)
		=
		\tr(G_Z|_{\tau^{-1}(x_0)},\Psi_Z|_{\tau^{-1}(x_0)}) \label{eq:cfsh}
	\end{align}
	
	Combining \eqref{eq:cfex} with \eqref{eq:cfsh}, we finally obtain the desired equality 
	\begin{align}
		\varphi^e_M(x_0)
		=
		\tr(G_Z|_{\tau^{-1}(x_0)},\Psi_Z|_{\tau^{-1}(x_0)})
		=
		\varphi^s_M(x_0).
	\end{align}
\end{proof}

Next, we shall describe the Lefschetz cycle $LC(G,\Psi)$ explicitly by using the constructible function $\varphi_M^s(G,\Psi)=\varphi_M^e(G,\Psi)$.

\begin{prp}[cf.\ \text{\cite[Proposition~5.6]{IMT15}}]\label{prp:conic}
	Under the condition \eqref{cd:6-1}, one has an equality 
	\begin{align}\label{LEQ} 
		LC(G,\Psi)=CC(\varphi^s_M(G,\Psi))
	\end{align}
	as Lagrangian cycles in $T^*M$.
\end{prp}

\begin{proof}
	Let $\pi_M \colon T^*M \longrightarrow M$ be the projection. 
	Fix a point $\mathring{x} \in M$ and compare the both sides of \eqref{LEQ} on a neighborhood of $\pi_M^{-1}(\mathring{x}) \subset T^*M$. 
	By the homotopy invariance of Lefschetz cycles (see Proposition~\ref{prp:ht}), \cite[Proposition~9.6.8]{K-S} and \eqref{ADDEE}, taking a sufficiently small open neighborhood $U$ of $\mathring{x}$ and replacing $\psi$ with $t\psi$ for $|1-t| \ll 1$, we may assume the following conditions:
	\begin{itemize}
		\item[(1)] $\cG|_U$ is trivial.
		\item[(2)] The open subset $U$ satisfies the condition of Proposition~\ref{prp:shrinking}.
		\item[(3)] $\Ev( \psi_x) \cap \{z \in \BC \mid |z|=1\}=\emptyset$ for any $x \in U$.
	\end{itemize}
	It suffices to show that 
	\begin{align}
		LC(G|_{\tl{U}},\Psi|_{\tl{U}}) 
		=
		CC(\varphi^s_M(G,\Psi)|_U),
	\end{align}
	where $\tl{U}:=\tau^{-1}(U)$. 
	As in the proof of Proposition~\ref{prp:compcf}, we can construct subbundles $\cG_+$ and $\cG_-$ of $\cG|_U$ for which we have the direct sum decomposition $\cG|_U =\cG_+ \oplus \cG_-$ and a metric on $\cG|_U$ satisfying \eqref{cd:metric}.
	By using this metric, we set 
	\begin{align}
		Z':=\{ (x,v_+,v_-) \in \cG|_U \mid |v_+| \le a, |v_-| < b \},
	\end{align}
	for some fixed constants $a,b>0$.
	Then $\psi^{-1}(Z) \cap Z$ is closed in $Z$ and open in $\psi^{-1}(Z)$ and hence we can construct a morphism
	\begin{align}
	\RG_{Z'}(\Psi) \colon \psi^{-1}(\RG_{Z'}(G)) \longrightarrow \RG_{Z'}(G) 
	\end{align}
	induced by $\Psi \colon \psi^{-1}G \longrightarrow G$. 
	Since $\cG_-$ is a shrinking subbundle of $\cG|_U$, by construction we have 
	\begin{align}
		\varphi^s_M(G,\Psi)|_U
		=
		\varphi(G_{\cG_-}^{-1!},\Psi_{\cG_-}^{-1!}). 
	\end{align}
	Moreover we have 
	\begin{align}
		CC(\varphi(G_{\cG_-}^{-1!},\Psi_{\cG_-}^{-1!}))
		=
		LC(G_{\cG_-}^{-1!},\Psi_{\cG_-}^{-1!}),
	\end{align}
	by Proposition~\ref{thm:5-8}.
	Thus it is enough to show that 
	\begin{align}
		LC(G|_{\tl{U}},\Psi|_{\tl{U}})
		=
		LC(G_{\cG_-}^{-1!},\Psi_{\cG_-}^{-1!})
	\end{align}
	as Lagrangian cycles in $T^*U$.
	In what follows, for simplicity we write $M$ instead of $U$ and $\cG|_U=\tl{U}, G|_{\tl{U}}, \Psi|_{\tl{U}}$ etc.\ by $\cG, G, \Psi$ etc., respectively.
	
	Let us take a $\mu$-stratification $\cG= \bigsqcup_{\alpha \in A} \cG_{\alpha}$ of $\cG$ which satisfies the following three conditions.
	\begin{enumerate}
		\item There exists a subset $B \subset A$ such that the zero-section $M \subset \cG$ of $\cG$ is $\bigsqcup_{\beta \in B} \cG_{\beta}$.
		\item $\MS(\RG_{Z'}(G)) \subset \bigsqcup_{\alpha \in A} T^*_{\cG_{\alpha}}\cG$ in $T^*\cG$.
		\item $\MS(G_{\cG_-}^{-1!}) \subset \bigsqcup_{\beta \in B}T^*_{\cG_{\beta}}M$ in $T^*M$.
	\end{enumerate}
	For $\beta \in B$, we denote $\cG_{\beta} \subset M$ by $M_{\beta}$. 
	Namely $M= \bigsqcup_{\beta \in B}M_{\beta}$ is a $\mu$-stratification of $M$. 
	Set $\Lambda= \bigsqcup_{\beta \in B} T^*_{M_{\beta}}M \subset T^*M$. 
	By the conditions above, we obtain
	\begin{align}
		\supp(LC(G,\Psi)), \ \supp(LC(G_{\cG_-}^{-1!},\Psi_{\cG_-}^{-1!})) \subset \Lambda.
	\end{align}
	Therefore it suffices to show that $LC(G,\Psi)$ coincides with $LC(G_{\cG_-}^{-1!},\Psi_{\cG_-}^{-1!})$ on an open dense subset of $\Lambda$. 
	Let $\Lambda_0$ be an open dense smooth subanalytic subset of $\Lambda$ whose decomposition $\Lambda_0= \bigsqcup_{i \in I}\Lambda_i$ into connected components satisfies the condition 
	\begin{align}
		\text{``for any $i \in I$, there exists $\beta_i \in B$ such that $\Lambda_i \subset T_{M_{\beta_i}}^*M$"}.
	\end{align}
	Fix $\Lambda_i$ and $M_{\beta_i}$ as above and we shall compare $LC(G,\Psi)$ with $LC(G_{\cG_-}^{-1!},\Psi_{\cG_-}^{-1!})$ on $\Lambda_i$.
	Take a point $p_0 \in \Lambda_i$ and set $x_0 :=\pi_M(p_0) \in M_{\beta_i}$. 
	Let $f \colon M \longrightarrow \BR$ be a real analytic function defined in an open neighborhood of $x_0$ which satisfies that $p_0=(x_0;df(x_0)) \in \Lambda_i$, $f(x_0)=0$ and the Hessian $\mathrm{Hess}(f|_{M_{\beta_i}})$ is positive definite. 
	Then by Corollary~\ref{crl:coeff}, we have 
	\begin{align}
		LC(G_{\cG_-}^{-1!},\Psi_{\cG_-}^{-1!})=m_i \cdot [T_{M_{\beta_i}}^*M]
	\end{align}
	in an open neighborhood of $\Lambda_i$ in $T^*M$, where $m_i \in \BC$ is defined by
	\begin{align}
		m_i
		:= 
		\dsum_{j \in \BZ} (-1)^j \tr\( H^j_{\{f \ge 0\}} (B(x_0,\delta);G_{\cG_-}^{-1!}) \xrightarrow{\Psi_{\cG_-}^{-1!}} H^j_{\{f \ge 0\}}(B(x_0,\delta);G_{\cG_-}^{-1!})\)
	\end{align}
	for a sufficiently small $\delta >0$. 
	Set $V:=B(x_0,\delta)$ and $W:=V \cap \{f<0\}$ in $M$. 
	Then 
	\begin{align}
		m_i
		=
		\tr(\RG_V(G_{\cG_-}^{-1!}),\RG_V(\Psi_{\cG_-}^{-1!}))
		-
		\tr(\RG_W(G_{\cG_-}^{-1!}),\RG_W(\Psi_{\cG_-}^{-1!})).\label{eq:6-16}
	\end{align}
	Set also $\tl{V}:=\tau^{-1}(V)$, $\tl{W}:=\tau^{-1}(W)\subset \cG$ and $\tl{f}:=f \circ \tau \colon \cG \longrightarrow \BR$.
	By Proposition~\ref{prp:shrinking} and the local invariance of characteristic classes, we obtain
	\begin{align*}
		\tr(\RG_V(G_{\cG_-}^{-1!}),\RG_V(\Psi_{\cG_-}^{-1!}))
		&=
		\int_{\tl{U}}C(\RG_{\tl V}(G),\RG_{\tl V}(\Psi)) \\
		&=
		\int_{\tl{U}} C(\RG_{\tl{V} \cap Z'}(G),\RG_{\tl{V}\cap Z'}(\Psi)) \\
		&=
		\tr(\RG_{\tl{V} \cap Z'}(G),\RG_{\tl{V}\cap Z'}(\Psi)).
	\end{align*}
	The last equality follows from the fact the support of $\RG_{\tl V \cap Z'}(G)$ is compact.
	Similarly, we get 
	\begin{align}
		\tr(\RG_W(G_{\cG_-}^{-1!}),\RG_W(\Psi_{\cG_-}^{-1!}))
		=
		\tr(\RG_{\tl{W} \cap Z'}(G),\RG_{\tl{W}\cap Z'}(\Psi)).
	\end{align}
	Applying Theorem~\ref{thm:idx} to the pair $(\RG_{\tl{V} \cap Z'}(G),\RG_{\tl{V} \cap Z'}(\Psi))$, we obtain
	\begin{align}\label{eq:6-23}
		\tr(\RG_V(G_{\cG_-}^{-1!}),\RG_V(\Psi_{\cG_-}^{-1!}))
		=
		\# ([\sigma_{f}] \cap LC(\RG_{\tl{V} \cap Z'}(G), \RG_{\tl{V} \cap Z'}(\Psi))).
	\end{align}
	Take a local coordinate $(x_1,\dots,x_m)$ in a neighborhood of $x_0$ such that $x_0=0$.
	Set $g(x):=|x|^2=x_1^2+\dots+x_m^2$.
	By the microlocal Bertini-Sard theorem (\cite[Proposition~8.3.12]{K-S}), there exists $\delta_0>0$ such that
	\begin{align}
		\Lambda \cap \Lambda_g \cap \pi_M^{-1}(\{0 <|x| \le \delta_0 \}) 
		=
		\emptyset, \label{cd:6-3} \\
		(\Lambda \widehat{+} T^*_{\{f=0\}}M) \cap \Lambda_g \cap \pi_M^{-1}(\{0 <|x| \le \delta_0 \}) 
		=
		\emptyset. \label{cd:6-4}
	\end{align}
	Define a real analytic function (defined on a neighborhood of $\tau^{-1}(x_0) \subset \cG$) $\tl{g}\colon \cG \longrightarrow \BR$ by $\tl{g}:=g \circ \tau$.
	Since $\tl{g}$ is proper on $\Supp(\RG_{\tl{V} \cap Z'}(G))$, by the microlocal Bertini-Sard theorem, there exists a sufficiently small $\delta_1>0$ such that
	\begin{align}
		\MS(\RG_{\tl{V} \cap Z'}(G)) \cap \Lambda_{\tl{g}} \cap \pi_{\cG}^{-1}(\{ v \in \cG \mid 0<|\tau(v)| \le \delta_1 \})
		=
		\emptyset, \label{cd:6-8} 
	\end{align}
	where $\pi_{\cG} \colon T^*\cG \longrightarrow \cG$ is the projection.
	Moreover by the proof of \cite[Theorem~9.5.6]{K-S}, there exists $\delta_2>0$ such that 
	\begin{align}
		c\ge0, 0<|x|\le \delta_2, f(x)>0 
		\Longrightarrow (x;c \cdot dg(x)+df(x)) 
		\not\in \Lambda. \label{cd:6-9}
	\end{align}
	Replacing the constant $\delta$ by a smaller one, we may assume that $0<\delta<\min(\delta_0,\delta_1,\delta_2)$.
	By the condition (i), \eqref{cd:6-8} and the definition of $\Lambda$, we have
	\begin{align}
		& 
		\lefteqn{
			\supp(LC(\RG_{\tl{V}\cap Z'}(G),\RG_{\tl{V}\cap Z'}(\Psi)))
		} \notag \\
		& \subset  
		\MS(\RG_{\tl{V} \cap Z'}(G)) \cap \cF_0 \label{eq:inc1} \\
		& \subset  
		\{ 
			\MS(\RG_{\tl V \cap Z'}(G)) \cup (\MS(\RG_{\tl V \cap Z'}(G)) +T^*_{\partial \tl{V}}\cG)
		\} \cap \cF_0 \\
		& \subset  
		\Lambda \cup (\Lambda + T^*_{\partial V}M)
		=:
		\Lambda^{\prime}. \label{eq:inc3}
	\end{align}
	Since $\Lambda^{\prime}$ is isotropic, by the microlocal Bertini-Sard theorem, there exists a sufficiently small $\e_0 >0$ such that 
	\begin{align}
		\Lambda^{\prime} \cap \Lambda_{f} \cap \pi_M^{-1}(\{0 < | f |\le \e_0 \})=\emptyset. \label{cd:6-10}
	\end{align}
	Arguing as in the proof of \cite[Theorem~9.5.6]{K-S} by using the conditions \eqref{cd:6-4}, \eqref{cd:6-9}, and \eqref{cd:6-10} and the estimates \eqref{eq:inc1}--\eqref{eq:inc3}, we obtain
	\begin{align}
		\Lambda_{f} \cap \supp(LC(\RG_{\tl{V} \cap Z'}(G), \RG_{\tl{V} \cap Z'}(\Psi))) 
		\subset 
		\pi_M^{-1}(\{f<-\e_0 \}) \sqcup \{p_0\}.
	\end{align}
	Hence from \eqref{eq:6-23} we deduce
	\begin{align}
		& 
		\lefteqn{
			\tr(\RG_V(G_{\cG_-}^{-1!})),\RG_V(\Psi_{\cG_-}^{-1!})))
		}  \notag \\
		& = 
		\# 
		\{
			\pi_M^{-1} (\{f <-\e_0 \}) \cap [\sigma_{f}] \cap LC(\RG_{\tl{V} \cap Z'}(G), \RG_{\tl{V} \cap Z'}(\Psi))
		\} \label{eq:6-38} \\
		& +
		[\sigma_{f}] \underset{p_0}{\cdot} LC(\RG_{Z'}(G), \RG_{Z'}(\Psi)),\notag
	\end{align}
	where the symbol $[\sigma_{f}] \underset{p_0}{\cdot} LC(\RG_{Z'}(G), \RG_{Z'}(\Psi))$ is the local intersection number of $[\sigma_{f}]$ and $LC(\RG_{Z'}(G), \RG_{Z'}(\Psi))$ at $p_0 \in \Lambda_i$. 

	The other term $\tr(\RG_W(G_{\cG_-}^{-1!}),\RG_W(\Psi_{\cG_-}^{-1!}))=\tr(\RG_{\tl{W} \cap Z'}(G),\RG_{\tl{W}\cap Z'}(\Psi))$ can be calculated as follows. 
	For $\e>0$, set $W_{\e}:=W \cap \{ f <-\e\}$ and $\tl{W_{\e}}:=\tl{W} \cap \{ \tl{f}<-\e\}=\tau^{-1}(W_{\e})$.

	\begin{lmm}\label{lem:6-3}
		There exists a sufficiently small $\e_1>0$ such that
		\begin{align}
			\tr(\RG_{\tl{W}\cap Z'}(G),\RG_{\tl{W} \cap Z'}(\Psi))
			=
			\tr(\RG_{\tl{W_{\e}} \cap Z'}(G),\RG_{\tl{W_{\e}}\cap Z'}(\Psi))
		\end{align}
		for any $0<\e< \e_1$. 
	\end{lmm}
	
	\begin{proof}
		Set $\Sigma:=\MS(\RG_{\tl{V}\cap Z'}(G)) \subset T^*\cG$. 
		Then by the microlocal Bertini-Sard theorem, there exists $\e_1>0$ such that
		\begin{align}
			\Sigma \cap \Lambda_{\tl{f}} \cap \pi^{-1}(\{-\e_1 \le \tl{f} <0\})=\emptyset.
		\end{align}
		Hence by the microlocal Morse lemma (\cite[Corollary 5.4.19]{K-S}), for $0<\e<\e_1$ we obtain
		\begin{align}
			\RG(\{\tl{f}<0\};\RG_{\tl{V}\cap Z'}(G)) 
			\simto 
			\RG(\{\tl{f}<-\e\};\RG_{\tl{V} \cap Z'}(G)).
		\end{align}
	\end{proof}

	Let us continue the proof of Proposition~\ref{prp:conic}. 
	By Lemma~\ref{lem:6-3} and Theorem~\ref{thm:idx}, we obtain
	\begin{align}\label{eq:6-27}
		\tr(\RG_W(G_{\cG_-}^{-1!}),\RG_W(\Psi_{\cG_-}^{-1!}))
		=
		\# ([\sigma_{f}] \cap LC(\RG_{\tl{W_{\e}}\cap Z'}(G),\RG_{\tl{W_{\e}}\cap Z'}(\Psi)))
	\end{align}
	for $0<\e<\e_1$. 
	Moreover it follows from the condition~(i) and the definition of $\Lambda$ that 
	\begin{align*}
		\supp(LC(\RG_{\tl{W_{\e}}\cap Z'}(G),\RG_{\tl{W_{\e}}\cap Z'}(\Psi)))
		& \subset 
		\MS(\RG_{\{\tl{f}<-\e\}}(\RG_{\tl{V}\cap Z'}(G))) \cap \cF_0 \\
		&\subset 
		\Lambda^{\prime} +\BR_{\leq 0}\Lambda_{f}.
	\end{align*}
	Comparing this last estimate with \eqref{cd:6-10}, we obtain
	\begin{align}
		\Lambda_{f} \cap \supp(LC(\RG_{\tl{W_{\e}}\cap Z'}(G),\RG_{\tl{W_{\e}}\cap Z'}(\Psi)))
		\subset 
		\pi_M^{-1}(\{f<-\e_0 \})
	\end{align}
	for $0<\e<\min(\e_0,\e_1)$. 
	Since
	\begin{align}
		LC(\RG_{\tl{W_{\e}}\cap Z'}(G), \RG_{\tl{W_{\e}}\cap Z'}(\Psi))
		=
		LC(\RG_{\tl{V}\cap Z'}(G),\RG_{\tl{V}\cap Z'}(\Psi))
	\end{align}
	on $\pi_M^{-1}(\{f<-\e_0 \})$, from \eqref{eq:6-27} we obtain 
	\begin{align}\label{eq:6-36}
		\lefteqn{
			\tr(\RG_W(G_{\cG_-}^{-1!}),\RG_W(\Psi_{\cG_-}^{-1!}))
		} \notag \\
		&=\#
		\{
			\pi_M^{-1}(\{f <-\e_0 \})\cap [\sigma_{f}] \cap LC(\RG_{\tl{V}\cap Z'}(G),\RG_{\tl{V}\cap Z'}(\Psi))
		\}.
	\end{align}
	Putting \eqref{eq:6-38} and \eqref{eq:6-36} into \eqref{eq:6-16}, we finally obtain
	\begin{align}
		m_i 
		=
		[\sigma_{f}] \underset{p_0}{\cdot} LC(\RG_{Z'}(G), \RG_{Z'}(\Psi)),
	\end{align}
	which shows
	\begin{align}
		LC(\RG_{Z'}(G), \RG_{Z'}(\Psi))
		=
		LC(G_{\cG_-}^{-1!},\Psi_{\cG_-}^{-1!}) \label{eq:6-40}
	\end{align}
	on $\Lambda_i$. 
	By the local invariance of Lefschetz cycles, we have 
	\begin{align}
		LC(\RG_{Z'}(G), \RG_{Z'}(\Psi))
		=
		LC(G,\Psi).\label{eq:6-41}
	\end{align}
	The result follows from \eqref{eq:6-40} and \eqref{eq:6-41}.
\end{proof}

We return to the original situation.
Namely, for a smooth fixed point component $M$ of $\phi$, we assume the condition
\begin{align}
	\text{``$1\notin \Ev(\phi^{\prime}_x)$ for any $x \in M$"}.
\end{align}
Setting $\cG=T_MX$ and $\psi=\phi'$, we can define the following constructible function by Proposition~\ref{prp:compcf}.

\begin{dfn}[\text{\cite[Definition 5.9]{IMT15}}]
	One defines a $\BC$-valued constructible function $\theta (F,\Phi)_M \in \CF(M)_{\BC}$ on $M$ by 
	\begin{align}
		\theta (F,\Phi)_M &:= \varphi^s_M(\nu_M(F),\Phi') \
		(=\varphi^e_M(\nu_M(F),\Phi')). 
	\end{align}
	and calls it the \textit{local trace function} of $(F, \Phi )$ on the fixed point component $M$. 
\end{dfn}

Combining Proposition~\ref{prp:specialization} with Proposition~\ref{prp:conic}, we have the following explicit description of the Lefschetz cycle $LC(F,\Phi)$.

\begin{thm}[\text{\cite[Theorem~5.10]{IMT15}}]\label{thm:main}
	In the situation as above, one has an equality
	\begin{align}
	LC(F,\Phi)
	=
	CC(\theta(F,\Phi)_M)
	\end{align}
	as Lagrangian cycles in $T^*M$.
\end{thm}

\begin{crl}\label{crl:cal}
	In the situation as above, assume moreover that $\Supp(F) \cap M$ is compact.
	Then the local contribution $c(F,\Phi)_M$ from $M$ can be calculated as the topological (or Euler) integral of the local trace function $\theta(F,\Phi)_M$:
	\begin{align}
	c(F,\Phi)_M
	=
	\int_M \theta(F,\Phi)_M.
	\end{align}
\end{crl}

\begin{proof}
	The result follows from the microlocal index formula (Theorem~\ref{thm:idx}) and Theorem~\ref{thm:main}.
\end{proof}

\section{Some examples}
In this section, we give some examples where local contributions can be explicitly determined by our method.

\begin{exa}[\text{cf.\ \cite[Example~6.5]{IMT15}}]\label{exa:proj}
	Let $X=\BR \BP^2$, $\phi \colon X \longrightarrow X, \ [x:y:z] \longmapsto [x:2y:z]$, and set 
	\begin{align}
		Z:=\{[x:y:z] \in \BR \BP^2 \mid x y z=0\} \subset \BR \BP^2.
	\end{align}
	Then $\phi$ induces a homeomorphism of $Z$.
	Hence, for the constructible object $F=\BC_Z$ on $\BR \BP^2$, we have a natural isomorphism $\Phi \colon \phi^{-1}F \simto F$.
	Let us consider the global Lefschetz number of the pair $(F,\Phi)$.
	Since
	\begin{align*}
		H^j(X;F)
		\simeq 
		H^j(Z;\BC)
		\simeq 
		\begin{cases}
		\BC & (j=0) \\
		\BC^4 & (j=1) \\
		0 & (\text{otherwise})
		\end{cases}
	\end{align*}
	and $\phi|_Z$ is homotopic to $\id_Z$, we can compute the global Lefschetz number directly as $\tr(F,\Phi)=1-4=-3$.
	
	Next, we shall calculate the local contributions. 
	The fixed point set $M$ of $\phi$ is decomposed into connected components as 
	\begin{align}
		M=\{y=0\} \sqcup \{[0:1:0]\}
		\simeq 
		\BR \BP^1 \sqcup \pt.
	\end{align}
	Let us first compute the local trace function $\theta(F,\Phi)_{\BR \BP^1}$ on $\{y=0\} \simeq \BR \BP^1$.
	Locally on $T_{\{y=0\}}X$, we can identify $(\nu_{\BR \BP^1}(F),\Phi')$ with $(F,\Phi)$.
	The differential $\phi'_x$ at any point $x \in \{y=0\} \simeq \BR \BP^1$ is equal to $2$.
	Hence we can take the whole normal bundle $T_{\{y=0\}}\BR \BP^2$ as an expanding subbundle and the zero-section as a shrinking subbundle (globally on $\{y=0\} \simeq \BR \BP^1$).
	For $x \neq [1:0:0],[0:0:1]$, the value is calculated as 
	\begin{align*}
		\theta(F,\Phi)_{\BR \BP^1}(x)
		& =
		\tr\left(\RG_c(\BR;\BC_0) \xrightarrow{(\phi_x)^*} \RG_c(\BR;\BC_0) \right) \\
		& =
		\tr\left( \BC \overset{\id}{\longrightarrow} \BC \right) 
		= 1.
	\end{align*}
	For $x =[1:0:0], [0:0:1]$, we have 
	\begin{align*}
		\theta(F,\Phi)_{\BR \BP^1}(x)
		& =
		\tr\left(
		\RG_c(\BR;\BC_\BR) \xrightarrow{(\phi_x)^*} \RG_c(\BR;\BC_\BR) \right) \\
		& =
		\tr\left(\BC[-1] \overset{\id}{\longrightarrow} \BC[-1]\right) 
		= -1.
	\end{align*}
	Therefore, the local trace function can be expressed as
	\begin{align}
		\theta(F,\Phi)_{\BR \BP^1}
		=
		1 \cdot \mathbf{1}_{\BR \setminus \{0\}} 
		+(-1) \cdot \mathbf{1}_{[1:0:0]}
		+(-1) \cdot \mathbf{1}_{[0:0:1]},
	\end{align}
	where $\mathbf{1}_Z$ denotes the characteristic function  of $Z$. 
	Thus by Corollary~\ref{crl:cal}, the local contribution from $\{y=0\} \simeq \BR \BP^1$ is computed as 
	\begin{align*}
		c(F,\Phi)_{\{y=0\}}
		& =
		\int_{\BR \BP^1} \theta(F,\Phi)_{\BR \BP^1} \\
		& =
		1 \cdot \chi_c(\BR \setminus \{0\})
		+(-1) \cdot \chi_c(\pt)
		+(-1) \cdot \chi_c(\pt) \\
		& =
		1 \cdot (-2)+(-1) \cdot 1+(-1) \cdot 1 \\
		& = -4,
	\end{align*}
	where $\chi_c(\cdot)$ denotes the Euler characteristic with compact support.
	Similar calculation shows that $c(F,\Phi)_{[0:1:0]}=1$.
	\medskip

	Next, we consider the self map $\psi \colon X \longrightarrow X, \ [x:y:z] \longmapsto [x:y/2:z]$, the same $\BR$-constructible sheaf $F=\BC_Z$ and the natural morphism $\Psi \colon \psi^{-1}F \simto F$.
	Also in this case, the global Lefschetz number $\tr(F,\Psi)$ is equal to $-3$.
	Moreover the fixed point set of $\psi$ is $\{y=0\}\sqcup \{[0:1:0]\} \simeq  \BR \BP^1 \sqcup \pt$ as in the case of $\phi$.
	The local trace functions on the connected components are the constant functions:
	\begin{align}
		\theta(F,\Phi)_{\BR\BP^1}
		=1 \cdot \mathbf{1}_{\BR \BP^1}, \quad
		\theta(F,\Phi)_{[0:1:0]}
		=-3 \cdot \mathbf{1}_{[0:1:0]}.
	\end{align}
	Thus we have 
	\begin{align}
		c(F,\Psi)_{\{y=0\}} =0, \quad
		c(F,\Psi)_{[0:1:0]} =-3.
	\end{align}
\end{exa}

\begin{exa}[\text{cf.\ \cite[Example~6.4]{IMT15}}]
	Let $S^2=\{x=(x_1, x_2, x_3) \in \BR^{3} \mid x_1^2+x_2^2 +x_{3}^2=1\}$ be the $2$-dimensional unit sphere in $\BR^{3}$ and $S^1= \{ e^{i \theta} \mid 0 \leq \theta \leq 2 \pi \}$ the $1$-dimensional one. 
	Set $X=S^1 \times S^2$. 
	For $e^{i \theta} \in S^1$ we define a real analytic isomorphism 
	$A_{\theta} \colon \BR^{3} \longrightarrow \BR^{3}$ of $\BR^{3}$ by
	\begin{align}
		A_{\theta} (x)=
		\begin{pmatrix}
			2 \cos \theta & -2 \sin \theta & 0\\
			2 \sin \theta & 2 \cos \theta & 0\\
			0 & 0 & 1
		\end{pmatrix}
		\begin{pmatrix}
		x_1\\x_2\\x_3
		\end{pmatrix}
	\end{align}
	and the one $\phi \colon X \longrightarrow X$ of $X$ by
	\begin{align}
		\phi(e^{i \theta}, x)
		=
		\left(e^{i \theta}, \dfrac{A_{\theta}(x)}{\|A_{\theta}(x)\|}\right).
	\end{align}
	Then the fixed point set $M$ of $\phi$ is a submanifold of $X$ and consists of three connected components $M_1,M_2, M_3$ defined by
	\begin{align}
		M_1=S^1 \times (0,0,1), \
		M_2=S^1 \times (0,0,-1), \
		M_3= \{ 1 \} \times (S^2 \cap \{ x_3=0 \}) \simeq S^{1},
	\end{align}
	respectively. 
	Note that for $p=( e^{i \theta}, (0,0,1)) \in M_1$ the set ${\rm Ev}( \phi^{\prime}_{p})$ of the eigenvalues of $\phi^{\prime}_{p} \colon (T_{M_1}X)_{p}  \longrightarrow (T_{M_1}X)_{p}$ is given by $\Ev( \phi^{\prime}_{p})= \{ 2e^{i \theta}, 2e^{-i \theta} \}$.
	In particular, it varies depending on the point $p \in M_1$ and satisfies the condition 
	\begin{equation}
		\text{``$1 \notin \Ev( \phi^{\prime}_{p})$ \quad for any $p \in M_1$"}.
	\end{equation}
	Let $\rho \colon S^2 \setminus \{ (0,0,1), (0,0,-1) \} \longrightarrow S^1$ be the natural surjective morphism and $\omega:=\exp(2\pi i/k) \in S^1$ a $k$-th root of unity.
	We denote the closure of $\rho^{-1}(\{1\} \sqcup \{\omega\} \sqcup \dots \sqcup \{\omega^{k-1}\})$ in $S^2$ by $K$. 
	Let us set
	\begin{align}
		Y= \left\{ (e^{i \theta},x) \in X \ \Big| \ x_3 > \frac{1}{2} \right\},
		\quad
		Z= \{ (1, x) \in Y \mid x \in K \}.
	\end{align}
	Then for the constructible sheaf $F=\BC_{Y \setminus Z} \in \Dc(X)$ there exists a natural morphism $\Phi : \phi^{-1}F \longrightarrow F$ and the inclusion map $i_{M_1} \colon M_1 \longhookrightarrow X$ is characteristic.
	Therefore in this case, we cannot apply \cite[Corollary 6.5]{M-T-3}.
	Since we can take the zero-section as a shrinking subbundle, we have
	\begin{align*}
		\theta(F,\Phi)_{M_1}((e^{i\theta},(0,0,1)))
		=
		\begin{cases}
			k & (\theta=0) \\
			1 & (\text{otherwise})
		\end{cases}
	\end{align*}
	for $0 \leq \theta < 2 \pi$. 
	By Corollary~\ref{crl:cal} we obtain
	\begin{align}
		c(F,\Phi) = c(F,\Phi)_{M_1}
		=k+ 1 \cdot (-1)= k-1.
	\end{align}
	On the other hand, we can easily see that
	\begin{align}
		\tr(F,\Phi) 
		=
		\chi_c(Y)-\chi_c(Z)
		=
		0-(1-k)=k-1.
	\end{align}
\end{exa}

\end{document}